\documentclass[12pt,leqno,a4paper,dvips,twoside]{article}

\usepackage{oldgerm}
\usepackage[french]{babel}
\usepackage{graphicx}
\usepackage{amsmath,amsfonts,amssymb}
\usepackage{makeidx}
\usepackage{shortvrb}
\usepackage{graphicx}
\usepackage{amsfonts}
\usepackage{amssymb}
\usepackage{amsmath}

\newtheorem{theo}{Th\'eor\`eme}%[section]
\newtheorem{cor}{Corollaire}%[section]
\newtheorem{prop}{Proposition}%[section]
\def\ds{\displaystyle}

\begin{document}

\begin{center}
%%% TITLE  %%%
{\large {\bf  Sur la suite des op\'erateurs Bernstein compos\'es }}
\vspace{0.5cm}

%%%% AUTOR %%%%
{\large Heiner Gonska (Duisburg-Essen) et Ioan Ra\c{s}a (Cluj-Napoca) }
\end{center}
\vspace{2.5cm}

{\small \noindent {\bf Abstrait:}
Nous consid\'erons une suite des op\'erateurs de Bernstein compos\'es et les formules de quadrature associ\'ees avec elles. Nous obtenons des bornes sup\'erieures pour l'erreur de l'approximation de fonctions continues et de l'approximation des integrales de fonctions continues. Les bornes sont donn\'ees en terme de modules de continuit\'e d'ordre un et deux. Deux in\'egalit\'es de type Tchebycheff-Gr\"uss sont aussi present\'ees.

\bigskip 

{\small \noindent {\bf }
\bigskip \bigskip


\noindent {\bf Mots cl\'es:} op\'erateurs de Bernstein compos\'es, formules de quadrature compos\'ees, modules de continuit\'e, degr\'e d'approximation, in\'egalit\'e de type Tche\-by\-cheff-Gr\"uss. 
}

\section{Introduction}
Dans l'article \cite{[BarMic]} D. B\v{a}rbosu et D. Micl\v{a}u\c{s} ont consider\'e une formule de quadrature bas\'ee sur des polyn\^omes de type Bernstein compos\'es. Ils ont donn\'e une in\'egalit\'e pour le reste de la formule de quadrature pour des fonctions dans la classe $C^2 [0,1]$, l'\'espace de fonctions d\'efinie sur l'intervalle $[0,1]$ ayant deux d\'eriv\'ees continues. Dans cet article nous utilisons les op\'erateurs introduits par les auteurs cit\'es pour approcher toutes les fonctions dans la classe $C[0,1]$, et nous donnons une \'evalution de l'erreur en utilisant le deuxi\`eme module de continuit\'e. 

De plus, nous \'etudions les it\'erations d'ordre $r$ des op\'erateurs lorsque $r \to \infty$.

Pour la formule de quadrature de B\v{a}rbosu et Micl\v{a}u\c{s}, nous trouvons l'ordre de grandeur du reste pour toutes les fonctions de l'\'espace $C[0,1]$. Notre note contient aussi deux r\'esultats du type Tchebycheff-Gr\"uss concernant la non-multiplicativit\'e de l'op\'erateur et de la formule de quadrature.

\section{Definition des op\'erateurs $\overline{B}_{n,m}$}

Rappelons les faits suivants:
\begin{itemize}
\item[1.] Pour $a,b \in \mathbb{R}, a < b$ et $f \in \mathbb{R}^{[a,b]}$ le polyn\^ome de Bernstein de degr\'e $n \in \mathbb{N}$ associ\'e avec $f$ est donn\'e par$$
B_n^{[a,b]} (f;x) = \frac{1}{(b-a)^n} \cdot \sum_{k=0}^n {n \choose k} (x-a)^k (b-x)^{n-k} \cdot f \left( a+k \frac{b-a}{n}\right).
$$
\item[2.] Pour $g \in C^2 [a,b]$ on a
$$
g(x) - B_n^{[a,b]} (g;x) = - \frac{(x-a)(b-x)}{2n} \cdot g'' (\xi_x) , \xi_x \in (a,b).
$$
\end{itemize}

Nous allons \'etudier la m\'ethode d'approximation suivante pour les fonctions continues definies sur $[0,1]$:

On divise $[0,1]$ en sous-intervalles $\left[\frac{k-1}{m},\frac{k}{m}\right]$, $k = 1,\ldots , m\in \mathbb{N}$. Sur $\left[\frac{k-1}{m}, \frac{k}{m}\right]$ nous consid\'erons
$$
B_{n,k} (f;x) := B_n^{\left[\frac{k-1}{m},\frac{k}{m}\right]} (f;x) = m^n \cdot \sum_{i=0}^n {n \choose i} \left( x - \frac{k-1}{m} \right)^i \left( \frac{k}{m} - x\right)^{n-i} f\left( \frac{k - n+i}{m\cdot n} \right).
$$
Maintenant nous composons les $B_{n,k} (f;\cdot)$ pour obtenir l'op\'erateur $\overline{B}_{n,m}$ defini par
$$
\overline{B}_{n,m} (f;x) := B_{n,k} (f;x) \;\;\mbox{si}\;\; x \in \left[\frac{k-1}{m} , \frac{k}{m}\right] , 1 \le k \le m .
$$
Ceci nous donne une fonction polyn\^omiale par morceaux de degr\'e $\le n$, continue aux points $\frac{k}{m} , 1 \le k \le m - 1$.

D'autre part, $\overline{B}_{n,m}$ est un operateur lin\'eaire et positif r\'eproduisant tous les fonctions lin\'eaires. Ces faits sont impliqu\'es par ceux de l'op\'erateur Bernstein classique (non-compos\'e).

\noindent Les op\'erateus $\overline{B}_{n,m}$ constituent une g\'en\'eralisation de
\begin{itemize}
\item l'op\'erateur de Bernstein sur $[0,1]$ - le cas $m = 1, n \in \mathbb{N}$,
\item l'interpolation lin\'eaire par morceaux $S_{\Delta_m}$ sur $[0,1]$ et aux points
$$
  \Delta_m : 0 < \frac{1}{m} < \frac{2}{m} < \ldots < \frac{m-1}{m} < 1.
$$
 - le cas $n = 1 , m \in \mathbb{N}$.
\end{itemize}

Chaque $\overline{B}_{n,m}$ est un cas sp\'ecial des operateurs spline de Schoenberg (''variation-diminishing spline operator'') associ\'es \`a une suite de n\oe uds appropri\'ee.

\section{Le degr\'e d'approximation par $\overline{B}_{n,m}$}

Plusieurs des nos resultats ci-dessous seront formul\'es \`a l'aide du module de continuit\'e d'ordre deux, donn\'e pour une fonction $f \in C[a,b]$ et $\delta \ge 0$ par
$$
\omega_2^{[a,b]}(f,\delta):= \sup\left\{ |f(x-h)-2f(x)+f(x+h)|: x+h, x-h \in [a,b], |h|\le \delta \right\}.
$$
Nous allons aussi utilizer la convention $\omega_{2} := \omega_{2}^{[0,1]}$.

Les deuxi\`emes moments $\overline{B}_{n,m} ((e_1-x)^2;x)$, o\`u $e_1(x) = x$, contr\^olent le degr\'e d'approximation. Pour $x \in \left[\frac{k-1}{m}, \frac{k}{m}\right]$ on a
$$
\overline{B}_{n,m} ((e_1-x)^2;x) = \frac{\left(x-\frac{k-1}{m}\right) \left(\frac{k}{m}-1\right)}{n}.
$$
Maintenant nous utilisons le r\'esultat suivant de P\v{a}lt\v{a}nea \cite{[Pal]}.

\begin{theo}

\noindent Si $L : C[0,1] \to C[0,1]$ est un operateur lin\'eaire et positif reproduisant toutes les fonctions lin\'eaires, alors pour tous $h > 0$ on a
$$
|L(f;x) - f(x)| \le \left[1 + \frac{1}{2h^2} \cdot L((e_1-x)^2;x)\right] \omega_2 (f;h).
$$
Si $L((e_1-x)^2;x) > 0$ le choix $h = \sqrt{L((e_1-x)^2;x)}$ implique
$$
|L(f;x) - f(x)| \le \frac32 \cdot \omega_2 (f;\sqrt{L((e_1-x)^2;x)});
$$
cette in\'egalit\'e est aussi valable si $L((e_1-x)^2;x) = 0$.
\end{theo}

Il en r\'esulte l'in\'egalit\'e suivante:
\begin{prop}
Pour $n,m \in \mathbb{N} , f \in C[0,1]$ et $x \in [0,1]$ on a
$$
|\overline{B}_{n,m} (f;x) - f(x)| \le \frac32 \omega_2 \left( f; \sqrt{\frac{\left(x - \frac{k-1}{m} \right) \left( \frac{k}{m} - x\right)}{n} } \right),
$$
pour $x \in \left[\frac{k-1}{m}, \frac{k}{m} \right]$, $1 \le k \le m$.
\end{prop}

\section{It\'erations des $\overline{B}_{n,m}$}

Consid\'erons $\overline{B}_{n,m}^\ell$, lorsque $\ell \to \infty$, avec $n,m$ fix\'es. 
Il est bien connu que chaque constituant
$$
B_{n,k} : C\left[\frac{k-1}{m},\frac{k}{m}\right] \to \Pi_n\big|_{\left[\frac{k-1}{m}, \frac{k}{m}\right]} , 1 \le k \le m ,
$$
produit une suite d'it\'erations $(B_{n,k})^\ell$, $\ell \ge 0$, qui pour toutes $f \in C \left[\frac{k-1}{m},\frac{k}{m}\right]$ g\'en\'ere une suite de polyn\^omes
$$
(B_{n,k})^\ell (f)
$$
qui approche, uniformement en $\left[\frac{k-1}{m},\frac{k}{m}\right]$, la fonction lin\'eaire $\ell_k$ interpolant $f$ aux points $\frac{k-1}{m}$ et $\frac{k}{m}$, c'est-\`a-dire, $\ell_k = B_{1,k} (f)$. D'ici il en r\'esulte que $(\overline{B}_{n,m})^\ell (f)$, $f \in C [0,1]$ converge uniformement vers $S_{\Delta_m}f$, l'interpolation lin\'eaire par morceaux.

En utilisant la transformation $\ell : [0,1] \to [a,b]$ donn\'ee par $\ell (x) = (b-a) x + a$, on peut \'ecrire
$$
\begin{array}{lcl}
B_n^{[a,b]} (f;x) & = & \ds \frac{1}{(b-a)^n} \sum_{k=0}^n {n \choose k} (x-a)^k (b-x)^{n-k} f\left( a+k \cdot \frac{b-a}{n}\right) \\[6mm]
& = & \ds B_n^{[0,1]} (f \circ \ell;y) \\[4mm]
& = & \ds \sum_{k=0}^n {n \choose k} y^k (1-y)^{n-k} (f\circ \ell) \left(\frac{k}{n}\right) \;\;\mbox{avec}\;\; y = \ell^{-1} (x) = \frac{x-a}{b-a} .
\end{array}
$$
Soit $r \in \mathbb{N}$, et consid\'erons l'it\'eration d'ordre $r$ de $B_n^{[a,b]}$, c'est \`a dire, $(B_n^{[a,b]})^r$. Nous utilisons le resultat suivant pour les it\'erations de $B_n = B_n^{[0,1]}$ donn\'e par Gonska, Kacs\'o et Pi\c tul dans l'article \cite{[GoKaPi]}.

\begin{prop}

Soit $B_n, n \in \mathbb{N}$, la suite des op\'erateurs de Bernstein classiques. Pour $r \in \mathbb{N}, \overline{f} \in C[0,1]$ et $x \in [0,1]$ on a
$$
|B_n^r (\overline{f};x) - B_1 (\overline{f};x)| \le \frac94 \cdot \omega_2 \left( \overline{f}, \sqrt{x(1-x) \left( 1 - \frac{1}{n}\right)^r}\right).
$$

\end{prop}

Ceci implique imm\'ediatement qu'on a, pour toutes $f \in C[a,b]$ et $x \in [a,b]$,
$$
\begin{array}{l}
\ds \left| \left(B_n^{[a,b]}\right)^r (f;x) - B_1^{[a,b]}(f;x)\right| \\[6mm]
\ds =\left|\left(B_n^{[0,1]}\right)^r \left(f\circ \ell ; \ell^{-1} (x)\right) - B_1^{[0,1]} (f\circ \ell , \ell^{-1}(x))\right| \\[6mm]
\ds \le \frac{9}{4} \cdot \omega_2^{[0,1]} \left( f \circ \ell , \sqrt{\ell^{-1} (x) \left[ 1-\ell^{-1}(x)\right] \left( 1 -\frac{1}{n}\right)^r}\right) \\[6mm]
\ds = \frac{9}{4} \cdot \omega_2^{[0,1]} \left( f \circ \ell , \sqrt{\frac{(x-a)}{(b-a)} \cdot \frac{(b-x)}{(b-a)} \cdot \left( 1 - \frac{1}{n} \right)^r}\right) \\[6mm]
\ds = \frac{9}{4} \cdot \omega_2^{[a,b]} \left( f; (b-a) \sqrt{\frac{(x-a)(b-x)}{(b-a)^2} \cdot \left( 1 - \frac{1}{n}\right)^r}\right) \\[6mm]
\ds = \frac{9}{4} \cdot \omega_2^{[a,b]} \left( f, \sqrt{(x-a)(b-x) \cdot \left( 1 - \frac{1}{n}\right)^r}\right) .
\end{array}
$$
Pour
$$
[a,b] = \left[ \frac{k-1}{m} , \frac{k}{m}\right] 
$$
et $\overline{f} : [0,1] \to \mathbb{R}$ consid\'erons la fonction
$$
f : = \overline{f}\big|_{[a,b]} : \left[ \frac{k-1}{m} , \frac{k}{m}\right] \to \mathbb{R} .
$$
Alors on en d\'eduit
$$
\begin{array}{l}
\ds \left|\left(\overline{B}_{n,m}\right)^r(\overline{f};x) - S_{\Delta_m} (\overline{f};x)\right| \\[6mm]
\ds \le \frac{9}{4} \omega_2^{\left[\frac{k-1}{m},\frac{k}{m}\right]} \left( \overline{f}\big|_{\left[\frac{k-1}{m},\frac{k}{m}\right]} , \sqrt{\left( x-\frac{k-1}{m}\right) \left(\frac{k}{m} - x\right) \left(1-\frac{1}{n}\right)^r}\right)  \\[6mm]
\ds = \frac{9}{4} \omega_2^{\left[\frac{k-1}{m},\frac{k}{m}\right]} \left( \overline{f} , \sqrt{\ldots}\right)  \\[6mm]
\ds \le \frac{9}{4} \omega_2^{[0,1]} \left( \overline{f} , \sqrt{\left( x - \frac{k-1}{m} \right)\left( \frac{k}{m} - x\right) \left( 1 - \frac{1}{n}\right)^r}\right) , \;\;\mbox{si}\;\; x \in \left[\frac{k-1}{m},\frac{k}{m}\right] , 1 \le k \le m . 
\end{array}
$$
Donc nous avons

\begin{prop}
Pour l'it\'eration d'ordre $r$ de l'op\'erateur $\overline{B}_{m,m} , \overline{f} \in C [0,1] , x \in [0,1]$ on a l'in\'egalit\'e
$$
\begin{array}{l}
\ds \left| \left(\overline{B}_{n,m}\right)^r (\overline{f},x) - S_{\Delta_m} (\overline{f},x)\right| \\[6mm]
\ds \le \frac{9}{4} \omega_2^{[0,1]} \left( \overline{f} , \sqrt{\left( x - \frac{k-1}{m} \right) \left( \frac{k}{m} - x \right) \left( 1 - \frac{1}{n} \right)^r}\right), \;\; x \in \left[\frac{k-1}{m} ,\frac{k}{m}\right] , 1 \le k \le m .
\end{array}
$$
Pour la norme uniforme il en r\'esulte
$$
\begin{array}{l}
\ds \left\| (\overline{B}_{n,m})^r (\overline{f}) - S_{\Delta_m} (\overline{f})\right\|_\infty \\[6mm]
\ds \le \frac{9}{4} \omega_2^{[0,1]} \left( \overline{f} ; \frac{1}{2m} \sqrt{\left( 1 - \frac{1}{n}\right)^r}\right) ,
\end{array}
$$
c'est \`a dire la convergence uniforme $(\overline{B}_{n,m})^r (\overline{f}) \to S_{\Delta_m} (\overline{f})$ pour $n,m$ fix\'es et $r \to \infty$.
\end{prop}

\section{Non-multiplicativit\'e de $\overline{B}_{n,m}$}

Dans cette section nous d\'emontrons une \'inegalit\'e de type Tchebycheff-Gr\"uss. Nous allons utiliser l'in\'egalit\'e generale suivante publi\'ee en \cite{[GoRaRu]}.

\begin{prop} Si $H : C[0,1] \to C[0,1]$ est un op\'erateur lin\'eaire et positif  reproduisant les fonctions constantes, alors pour toutes $f, g \in C[0,1]$ et $x \in [0,1]$ on a:
$$
\begin{array}{l}
T(f,g;x) := |H(f,g;x) - H(f;x) \cdot H(g;x)| \\[6mm]
\ds \le \frac14 \widetilde{\omega} \left( f; 2 \sqrt{H((e_1-x)^2;x)}\right) \widetilde{\omega} \left( g;2 \sqrt{H((e_1-x)^2;x)}  \right). \end{array} 
$$
Pour $t \in [0,\infty)$ la quantit\'e
$$
\omega(f;t) = \sup \{|f(x) - f(y)| : |x-y| \le t\}
$$
est le module de continuit\'e d'ordre un, et le plus petit majorant concave du module est donn\'e par $$
\widetilde{\omega} (f;t) = \left\{ 
\begin{array}{lll}

\ds \sup_{0 \le x \le t \le y \le 1 \atop x \not= y}{\frac{(t-x)\omega(f,y)+(y-t)\omega(f,x)}{y-x}}  & \mbox{ } & ,0 \le t \le 1, \\[6mm]
\omega (f,1) & \mbox{ } & ,t > 1 . 

\end{array} \right.
$$
\end{prop}

En substituant dans l'in\'egalit\'e la repr\'esentation des moments d'ordre deux de $\overline{B}_{n,m}$ on obtient
\begin{prop}
Pour $f,g \in C[0,1]$ et $x \in [0,1]$ l'in\'egalit\'e suivante de type Gr\"uss est valable:
$$
\begin{array}{l}
\left| \overline{B}_{n,m} (f \cdot g;x) - \overline{B}_{n,m} (f;x) \overline{B}_{n,m} (g;x)\right| \\[6mm]
\ds  \le \frac14 \widetilde{\omega} \left( f ; 2 \sqrt{\frac{\left( x - \frac{k-1}{m}\right) \left( \frac{k}{m} - x\right)}{n}} \right) \widetilde{\omega} \left( g; 2 \sqrt{\frac{\left(x - \frac{k-1}{m}\right) \left( \frac{k}{m} - x\right)}{n}} \right) , \;\; \mbox{si}\;\; x \in \left[ \frac{k-1}{m} , \frac{k}{m}\right] .
\end{array}
$$
\end{prop}

Remarquons que l' in\'egalit\'e au dessus refl\'ete le fait que $\overline{B}_{n,m}$ interpole aux points $\frac{k}{m}$, $0 \le k \le m$.

\section{Sur la formule de quadrature bas\'ee sur $\overline{B}_{n,m}$}

La formule de quadrature introduite par B\v{a}rbosu et Micl\v{a}u\c{s} est donn\'ee par
$$
\begin{array}{lcl}
\ds \int_0^1 f(x) dx = \sum_{k=1}^m \int_{\frac{k-1}{m}}^{\frac{k}{m}} f(x) dx & \approx & \ds \sum_{k=1}^m \int_{\frac{k-1}{m}}^{\frac{k}{m}} B_{n,k} (f;x) dx \\[6mm]
 & = & \ds \int_0^1 \overline{B}_{n,m} (f;x) dx =: I_{n,m} (f) .
 \end{array}
$$

\begin{theo}
Pour la formule de quadrature au-dessus on a
$$
I_{m,n} (f) = \frac{1}{m(n+1)} \sum_{k=1}^m \; \sum_{i=0}^n f \left( \frac{kn -  n + i}{m+n} \right).
$$
\end{theo}

\noindentÊ{\bf D\'emonstration:}
\noindent Soit $k$ fix\'e. On peut \'ecrire
$$
\begin{array}{l}
\ds \int_{\frac{k-1}{m}}^{\frac{k}{m}} B_{n,k} (f;x) dx \\[6mm]
\ds = m^n \sum_{i=0}^n {n \choose i} f \left( \frac{kn-n+i}{m\cdot n} \right) \int_{\frac{k-1}{m}}^{\frac{k}{m}} \left( x - \frac{k-1}{m} \right)^i \left( \frac{k}{m} - x \right)^{n-i} dx \\[6mm]
\ds = m^n \sum_{i=0}^n {n \choose i} f \left( \frac{kn-n+i}{m\cdot n} \right) \int_0^1 \left( \frac{i}{m} \right) \left[ \frac{1}{m} (1-t)\right]^{n-i} \frac{1}{m} dt \\[6mm]
\ds = \frac{1}{m} \sum_{i=0}^n {n\choose i} \int_0^1 t^i (1-t)^{n-i} dt \cdot f \left( \frac{kn - n + i}{m\cdot n} \right) \\[6mm]
\ds = \frac{1}{m} \sum_{i=0}^n {n \choose i} B (i+1,n-i+1) \cdot f\left( \frac{kn -n + i}{m \cdot n} \right) \\[6mm]
\ds = \frac{1}{m} \sum_{i=0}^n {n \choose i} \frac{1}{n+1} {n \choose i}^{-1} f \left( \frac{kn - n + i}{m \cdot n} \right) \\[6mm]
\ds = \frac{1}{m(n+1)} \sum_{i=0}^n f \left( \frac{kn-n+i}{m\cdot n} \right).
\end{array}
$$
Une sommation pour toutes les valeurs de $k$ donne la r\'epresentation desir\'ee. \hfill $\square$

\bigskip \bigskip

Le r\'esultat suivant est une am\'elioration significative du Theorem 2.2 de \cite{[BarMic]}. 

\begin{theo}
Pour $g \in C^2 [0,1]$ on a
$$
\left| \int_0^1 g(x) dx - I_{m,n} (g)\right| \le \frac{1}{12m^2n} \cdot \| g''\|_\infty .
$$
\end{theo}

\noindent {\bf D\'emonstration.} La preuve r\'esulte des (in)\'egalit\'es suivantes:
$$
\begin{array}{l}
\ds \left| \int_0^1 g(x) - I_{m,n} (g) \right| \\[6mm]
\ds = \Bigg| \sum_{k=1}^m \int_{\frac{k-1}{m}}^{\frac{k}{m}} g(x) dx - \sum_{k=1}^m \underbrace{\frac{1}{m(n+1)} \sum_{i=0}^n g \left( \frac{kn-n+i}{m\cdot n} \right)}_{= \ds \int_{\frac{k-1}{m}}^{\frac{k}{m}} B_{n,k} (g;x)dx} \\[6mm]
\ds = \left| \sum_{k=1}^m \int_{\frac{k-1}{m}}^{\frac{k}{m}} [g(x) - B_{n,k} (g;x)|dx \right| \\[6mm]
\ds = \left| \sum_{k=1}^m \int_{\frac{k-1}{m}}^{\frac{k}{m}} - \frac{\left( x - \frac{k-1}{m}\right) \left( \frac{k}{m} - y \right)}{2n} g'' (\xi_{x,k}) dx \right| , \xi_{x,k} \in \left(\frac{k-1}{m} , \frac{k}{m} \right) \\[6mm]
\ds \le \frac{1}{2n} ||g''||_\infty \sum_{k=1}^{m} \int_{\frac{k-1}{m}}^{\frac{k}{m}} \left( x - \frac{k-1}{m} \right) \left( \frac{k}{m} - x\right) dx  \\[6mm]
\ds = \frac{1}{2n} ||g''||_\infty \sum_{k=1}^m \frac{1}{m} \int_0^1 \frac{t}{m} \cdot \frac{1-t}{m} dt \\[6mm]
\ds = \frac{1}{2n \cdot m^2} ||g''||_\infty \cdot \int_0^1 t(1-t)dt \\[6mm]
\ds = \frac{1}{12nm^2} ||g''||_\infty .
\end{array} 
$$
\hfill $\square$
\bigskip

\noindent Dans le th\'eor\`eme suivant nous utilisons la fonctionnelle $K$ definie par
$$
K \left( \delta, f; C^0[0,1],C^2[0,1] \right) :=   \inf \left\{ ||f-g||_\infty + \delta||g''||_\infty : g \in C^2 [0,1]\right\}, \delta \ge 0\\[6mm].
$$

\begin{theo}
Pour $f \in C[0,1]$ et pour $m,n \ge 1$ on a
\begin{itemize}
\item[(i)] $\ds \left| \int_0^1 f(x) dx - I_{m,n} (f) \right| \le 2 K \left( \frac{1}{24m^2n},f; C^0[0,1],C^2[0,1]\right) $;
\item[(ii)] $\left|\int_0^1 f(x) dx - I_{m,n} (f) \right| \le \frac{9}{4} \omega_2 \left( f; \frac{1}{m\sqrt{6n}} \right)$.
\end{itemize}
\end{theo}

\bigskip

\noindent {\bf D\'emonstration.} Pour chaque $f \in C[0,1]$ on a 
$$
\left| \int_0^1 f(x) dx - I_{m,n} (f) \right| \le ||f||_\infty + \frac{1}{m(n+1)} \sum_{k=1}^m \sum_{i=0}^n ||f||_\infty = 2 ||f||_\infty .
$$
Alors, quelle que soit $g \in C^2 [0,1]$, nous deduisons, en notant $H(f) := \int_0^1 f(x) dx - I_{m,n} (f)$, que
$$
\begin{array}{lcl}
|H(f)| & = & |H(f-g+g)| \\[4mm]
& \le & |H(f-g)| + |H(g)| \\[2mm]
& \le & \ds 2 ||f-g||_\infty + \frac{1}{12nm^2} ||g''||_\infty .
\end{array}
$$
Il en r\'esulte
$$
\begin{array}{lcl}
|H(f)| & \le & \ds 2 \cdot \inf \left\{ ||f-g||_\infty + \frac{1}{24nm^2} ||g''||_\infty : g \in C^2 [0,1]\right\}\\[6mm]
& = & 2 \cdot K \left( \frac{1}{24m^2n} , f; C^0[0,1],C^2[0,1] \right) .
\end{array}
$$
\hfill$\square$
\vspace{10mm}

Pour d\'emontrer (ii) nous citons le Th\'eor\`eme 4.2 de Gonska et Kovacheva \cite{[GonKov]}.

\begin{theo} Soit $(B,||\cdot||_B)$ un \'espace de Banach , et soit $H : C[a,b] \to B$ un operator (pas n\'ecessairement lin\'eaire, pas n\'ecessairement positif) satisfaisant les conditions suivantes avec des constantes $\gamma , \alpha , \beta_0, \beta_1, \beta_2 {\ge 0}$ ind\'ependantes de $f$ et $g$:
\begin{enumerate}
\item[a)] $||H(f+g)||_B \le \gamma \{||Hf||_B + ||Hg||_B\}$ pour toute $f \in C[a,b]$,
\item[b)] $||Hf||_B \le \alpha ||f||_\infty$ pour toute $f \in C[a,b];$
\item[c)] $||Hg||_B \le \beta_0 ||g||_\infty + \beta_1 ||g'||_\infty + \beta_2 ||g''||_\infty$ pour toute $g \in C^2 [a,b]$.
\end{enumerate}

\bigskip

Alors, quelque soit $f \in C[a,b] , 0 < h \le \frac{b-a}{2}$, nous avons
$$
||Hf||_B \le \gamma \left\{ \beta_0 ||f||_\infty + \frac{2\beta_1}{h} \omega_1 (f;h) + \frac34 \left( \alpha + \beta_0 + \frac{2 \beta_1}{h} + \frac{2 \beta_2}{h^2} \right) \omega_2 (f;h)\right\}.
$$
\end{theo}

Dans le cas pr\'esent nous prenons
$$
\begin{array}{l}
C[a,b] = C[0,1] , B = \mathbb{R}, \\[4mm]
\gamma = 1 , \alpha = 2 , \beta_0 = 0 , \beta_1 = 0 , \beta_2 = \frac{1}{12m^2n} .
\end{array}
$$
On obtient, pour $0 < h \le \frac12$,
$$
|H(f)| \le \frac34 \left( 2 + \frac{1}{h^2} \cdot \frac{1}{6m^2n} \right) \omega_2 (f,h).
$$
En choisissant $h = \frac{1}{\sqrt{6m^2n}}$ nous arrivons \`a (ii). \hfill $\square$

%%%
\section{Non-multiplicativit\'e de la formule de quadrature}

Considerons maintenant de nouveau la formule de quadrature
$$
I_{n,m} (f) = \int_0^1 \overline{B}_{n,m} (f;x) dx.
$$
Ici, notre but est de donner une borne sup\'erieure pour la quantit\'e
$$
|T(f,g)| := |I_{n,m} (f\cdot g) - I_{n,m} (f) I_{n,m} (g)|.
$$
\`A cette fin, nous utilisons de nouveau le majorant concave $\widetilde{\omega}$ et  le resultat suivant de \cite{[GoRaRu]}, Th. 3.1. 

\begin{prop}
Si $L: C[0,1] \to \mathbb{R}$ est une fonctionnelle lin\'eaire et positive satisfaisant $L(e_0) = 1$, alors pour toutes $f,g \in C[0,1]$ nous avons$$
|T(f,g)| \le \frac14 \widetilde{\omega} (f; 2 \sqrt{T(e_1,e_1)}) \widetilde{\omega} (g; 2 \sqrt{T(e_1,e_1)}).
$$
Ici,
$$
T(e_1,e_1) = L(e_2) - [L(e_1)]^2.
$$
\end{prop}

\bigskip

La proposition ci-dessus conduit \`a
\begin{prop}
$$
\begin{array}{l}
\ds |I_{n,m}(f\cdot g) - I_{n,m} (f) I_{n,m} (g) | \\[6mm]
\ds \le \frac14 \widetilde{\omega} \left( f; 2 \sqrt{\frac{1}{12} + \frac{1}{6m^2n}} \right) \widetilde{\omega} \left( g; 2 \sqrt{\frac{1}{12} + \frac{1}{6m^2 n}}\right) .
\end{array}
$$
\end{prop}

\noindent {\bf D\'emonstration.} Il suffit de calculer
$$
\begin{array}{lcl}
T(e_1,e_1) & = & \ds I_{n,m} (e_2) - [I_{n,m} (e_1)]^2 \\[4mm]
& = & \ds \int_0^1 \overline{B}_{n,m} (e_2;x) dx - \left[ \int_0^1 \overline{B}_{n,m} (e_1;x) dx \right]^2 \\[4mm]
& = & \ds \int_0^1 \overline{B}_{n,m} (e_2;x) dx - \frac14 \\[4mm]
& = & \ds - \frac14 + \sum_{k=1}^m \int_{\frac{k-1}{m}}^{\frac{k}{m}} B_{n,k} (e_2;x)dx \\[4mm]
& = & \ds - \frac14 + \frac{1}{m(n+1)} \sum_{k=1}^m \sum_{i=0}^n \left( \frac{kn - n + i}{m\cdot n}\right)^2 \\[4mm]
& = & \ds - \frac14 + \frac{1}{m^3n^2(n+1)} \sum_{k=1}^m \sum_{i=0}^n (kn-n+i)^2 \\[6mm]
& = & \ds \frac{1}{12} + \frac{1}{6m^2n} .
\end{array}
$$
\vspace{-20mm}
\mbox{} \hfill $\square$

\vspace{20mm}

\bigskip

\begin{cor}
Pour $n,m \to \infty$ nous obtenons
$$
\begin{array}{l}
\ds \left| \int_0^1 (f\cdot g)(x) dx - \int_0^1 f(x) dx \int_0^1 g(x)dx \right| \\[4mm]
\ds = \left| \lim_{n,m\to \infty} \left\{ \int_0^1 \overline{B}_{n,m} (f\cdot g;x) dx - \int_0^1 \overline{B}_{n,m} (f;x) dx \int_0^1 \overline{B}_{n,m} (g;x)dx \right\}Ê\right| \\[4mm]
\ds \le \lim_{n,m \to \infty} \frac14 \widetilde{\omega} (f; 2 \sqrt{\frac{1}{12} + \frac{1}{6m^2n} }) \widetilde{\omega} (g;2\sqrt{\frac{1}{12} + \frac{1}{6m^2n} }) \\[4mm]
\ds =\frac{1}{4} \widetilde{\omega} \left( f; \frac{1}{\sqrt{3}}\right) \widetilde{\omega} \left( g; \frac{1}{\sqrt{3}} \right) \\[4mm]
\ds \le \frac{1}{12} ||f'||_{L_\infty} ||g'||_{L_\infty} \;\; \mbox{pour} \;\; f',g' \in L_\infty [0,1].
\end{array}
\vspace{-10mm}
$$\hfill $\square$

\vspace{10mm}

Remarquons qu'il s'agit d'une in\'egalit\'e de type Tchebycheff-Gr\"uss pour la fonctionnelle d'int\'egration o\`u la constante $\frac{1}{12}$ est la meilleure possible. 
\end{cor}

\noindent
{\bf Remerciement: }Les auteurs remercient chaleuresement Mme Birgit Dunkel pour la r\'ealisation de ce manuscrit. On remercie aussi M. Catalin Badea pour quelques corr\'ections linguistiques.

\bigskip

\noindent
 $\begin{array}{ll}
\textrm{Heiner Gonska}\\
 \textrm{University of Duisburg-Essen} \\
 \textrm{Department of Mathematics} \\
 \textrm{D-47048 Duisburg} \\
 \textrm{Germany}\\
\textrm{e-mail: heiner.gonska@uni-due.de} \end{array}
$

\bigskip
 \noindent
 $\begin{array}{ll}
\textrm{Ioan Ra\c{s}a}\\
 \textrm{Technical University} \\
 \textrm{Department of Mathematics } \\
 \textrm{RO-400020 Cluj-Napoca} \\
 \textrm{Romania}\\
\textrm{e-mail: Ioan.Rasa@math.utcluj.ro} \end{array} $

\bigskip
\noindent

\end{document}